# CONSISTENCY OF CROSS VALIDATION FOR COMPARING REGRESSION PROCEDURES[1]

By Yuhong Yang

*University of Minnesota*

Theoretical developments on cross validation (CV) have mainly focused on selecting one among a list of finite-dimensional models (e.g., subset or order selection in linear regression) or selecting a smoothing parameter (e.g., bandwidth for kernel smoothing). However, little is known about consistency of cross validation when applied to compare between parametric and nonparametric methods or within nonparametric methods. We show that under some conditions, with an appropriate choice of data splitting ratio, cross validation is consistent in the sense of selecting the better procedure with probability approaching 1.

Our results reveal interesting behavior of cross validation. When comparing two models (procedures) converging at the same nonparametric rate, in contrast to the parametric case, it turns out that the proportion of data used for evaluation in CV does not need to be dominating in size. Furthermore, it can even be of a smaller order than the proportion for estimation while not affecting the consistency property.

**1. Introduction.** Cross validation (e.g., Allen [2], Stone [25] and Geisser [9]) is one of the most commonly used model selection criteria. Basically, based on a data splitting, part of the data is used for fitting each competing model (or procedure) and the rest of the data is used to measure the performance of the models, and the model with the best overall performance is selected. There are a few different versions of cross-validation (CV) methods, including delete-1 CV, delete-$k$ ($k > 1$) CV and also generalized CV methods (e.g., Craven and Wahba [6]).

Cross validation can be applied to various settings, including parametric and nonparametric regression. There can be different primary goals when

Received March 2006; revised February 2007.
[1]Supported in part by NSF Career Grant 0094323.

applying a CV method: one mainly for identifying the best model/procedure among the candidates and another mainly for estimating the mean function or for prediction (see, e.g., Geisser [9]). A number of theoretical results have been obtained, mostly in the areas of linear regression and in smoothing parameter selection for nonparametric regression. In linear regression, it has been shown that delete-1 and generalized CVs are asymptotically equivalent to the Akaike Information Criterion (AIC) [1] and they are all inconsistent in the sense that the probability of selecting the true model does not converge to 1 as $n$ goes to $\infty$ (see Li [15]). In addition, interestingly, the analysis of Shao [19] showed that in order for delete-$k$ CV to be consistent, $k$ needs to be dominatingly large in the sense that $k/n \to 1$ (and $n-k \to \infty$). Zhang [35] proved that delete-$k$ CV is asymptotically equivalent to the Final Prediction Error (FPE) criterion when $k \to \infty$. The readers are referred to Shao [20] for more asymptotic results and references on model selection for linear regression. In the context of nonparametric regression, delete-1 CV for smoothing parameter selection leads to consistent regression estimators (e.g., Wong [32] for kernel regression and Li [14] for the nearest-neighbor method) and leads to asymptotically optimal or rate-optimal choice of smoothing parameters and/or optimal regression estimation (see, e.g., Speckman [22] and Burman [5] for spline estimation, Härdle, Hall and Marron [12], Hall and Johnstone [11] and references therein for kernel estimation). Györfi et al. [10] gave risk bounds for kernel and nearest-neighbor regression with bandwidth or neighbor size selected by delete-1 CV. See Opsomer, Wang and Yang [17] for a review and references related to the use of CV for bandwidth selection for nonparametric regression with dependent errors.

In real-world applications of regression, in pursuing a better estimation accuracy, one may naturally consider the use of cross validation to choose between a parametric estimator and a nonparametric estimator (or at least to understand their relative performance). Similarly, when different types of nonparametric estimators are entertained as plausible candidates, cross validation is also applicable to choose one of them. Recently, a general CV methodology has been advocated by van der Laan, Dudoit, van der Vaart and their co-authors (e.g., van der Laan and Dudoit [26], van der Laan, Dudoit and van der Vaart [27] and van der Vaart, Dudoit and van der Laan [28]), which can be applied in other contexts (e.g., survival function estimation). Risk bounds for estimating the target function were derived and their implications on adaptive estimation and asymptotic optimality were obtained. When CV is used for the complementary purpose of identifying the best candidate, however, it is still unclear whether CV is generally consistent and if the data splitting ratio has a sensitive effect on consistency. For successful applications of CV in practice, a theoretical understanding of these issues is very much of interest.

In this paper, we address the aforementioned consistency issue and show that a voting-based cross validation is consistent for comparing general regression procedures when the data splitting ratio is properly chosen.

In the context of linear regression, Shao's result [19] implies that with $k/n$ not converging to 1, delete-$k$ CV does not differentiate well between two correct models. However, in nonparametric regression, it turns out that this is not the case. In fact, as long as at least one of the competing procedures converges at a nonparametric rate, the estimation size and evaluation size can be of the same order in data splitting, and sometimes the estimation size can even be the dominating one.

In the settings where theoretical properties of CV were investigated before, the best model (in the linear regression context) or the best smoothing parameter (in the case of kernel regression) exists in a natural way. When comparing two general estimators, the issue becomes more complicated. In this paper, we compare two estimators in terms of a loss function and the consistency in selection is established when one estimator is better than the other in that sense.

The paper is organized as follows. In Section 2 we set up the problem. The main result is presented in Section 3, followed by simulation results in Section 4. Concluding remarks are in Section 5. The proof of the main result is in Section 6.

## 2. Problem setup.

Consider the regression setting

$$Y_i = f(X_i) + \varepsilon_i, \qquad 1 \le i \le n,$$

where $(X_i, Y_i)_{i=1}^n$ are independent observations with $X_i$ i.i.d. taking values in a $d$-dimensional Borel set $\mathcal{X} \subset R^d$ for some $d \ge 1$, $f$ is the true regression function and $\varepsilon_i$ are the random errors with $E(\varepsilon_i|X_i) = 0$ and $E(\varepsilon_i^2|X_i) < \infty$ almost surely. The distribution of $X_i$ is unknown.

Rates of convergence of various popular regression procedures have been well studied. Under the squared $L_2$ loss (and other closely related performance measures as well), for parametric regression (assuming that the true regression function has a known parametric form), estimators based on maximum likelihood (if available) or least squares usually converge at the rate $n^{-1}$. For nonparametric estimation, the convergence rates are slower than $n^{-1}$ with the actual convergence rate depending on both the regression procedure and the smoothness property of the true regression function.

In applications, with many regression procedures available to be applied, it is often challenging to find the one with best accuracy. One often faces the issue: Should I go parametric or nonparametric? Which nonparametric procedure to use? As is well known, nonparametric procedures have more flexibility yet converge suboptimally compared to estimation based on a correctly specified parametric model. Sometimes, for instance, there is a clear linear trend and naturally a simple linear model is a good candidate. It may be unclear, however, whether a nonparametric method can provide a better estimate to capture a questionable slight curvature seen in the data. For a specific example, for the rather famous Old Faithful Geyser data (Weisberg [31]), there were several analyses related to the comparison of linear regression with nonparametric alternatives (see, e.g., Simonoff [21] and Hart [13]).

For simplicity, suppose that there are two regression procedures, say $\delta_1$ and $\delta_2$, that are considered. For example, $\delta_1$ may be simple linear regression and $\delta_2$ may be a local polynomial regression procedure (see, e.g., Fan and Gijbels [8]). For another example, $\delta_1$ may be a spline estimation procedure (see, e.g., Wahba [29]) and $\delta_2$ may be a wavelet estimation procedure (see, e.g., Donoho and Johnstone [7]). Based on a sample $(X_i, Y_i)_{i=1}^n$, the regression procedures $\delta_1$ and $\delta_2$ yield estimators $\widehat{f}_{n,1}(x)$ and $\widehat{f}_{n,2}(x)$, respectively. We need to select the better one of them.

Though for simplicity we assumed that there are only two competing regression procedures, similar results hold when a finite number of candidate regression procedures are in competition. We emphasize that our focus in this work is not on tuning a smoothing parameter of a regression procedure (such as the bandwidth for kernel regression). Rather, our result is general and applicable also for the case when the candidate regression procedures are very distinct, possibly with different rates of convergence for estimating the regression function.

Cross validation is a natural approach to address the above model/procedure comparison issue. It has the advantage that it requires mild distributional assumptions on the data and it does not need to find characteristics such as degrees of freedom or model dimension for each model/procedure (which is not necessarily easy for complicated adaptive nonparametric procedures or even parametric procedures where model selection is conducted). To proceed, we split the data into two parts: the estimation data consist of $Z^1 = (X_i, Y_i)_{i=1}^{n_1}$ and the validation data consist of $Z^2 = (X_i, Y_i)_{i=n_1+1}^{n}$. Let $n_2 = n - n_1$. We apply $\delta_1$ and $\delta_2$ on $Z^1$ to obtain the estimators $\widehat{f}_{n_1,1}(x)$ and $\widehat{f}_{n_1,2}(x)$, respectively. Then we compute the prediction squared errors of the two estimators on $Z^2$:

$$\text{CV}(\widehat{f}_{n_1,j}) = \sum_{i=n_1+1}^{n} (Y_i - \widehat{f}_{n_1,j}(X_i))^2, \qquad j = 1, 2. \tag{2.1}$$

If $\text{CV}(\widehat{f}_{n_1,1}) \le \text{CV}(\widehat{f}_{n_1,2})$, $\delta_1$ is selected and otherwise $\delta_2$ is chosen. We call this a delete-$n_2$ CV. For investigating the consistency property of CV, in this work we only consider the case when $\min(n_1, n_2) \to \infty$.

A more balanced use of the data is CV with multiple data splittings. A voting-based version will be considered as well.

## 3. Consistency of cross validation.

In this section, we first define the concept of consistency when two general regression procedures are compared and then state the conditions that are needed for our result. Note that there are two types of consistency results on CV. One is from the selection perspective (which concerns us if the true model or the best candidate procedure is selected with probability tending to 1) and the other is in terms of convergence of the resulting estimator of the regression function (either in probability or in risk under a proper loss function). The former is our focus in this work. Note that in general, neither notion guarantees the other. Nevertheless, for a global loss function, consistency in selection usually implies consistency in estimation of the regression function at least in probability. In contrast, consistency in estimation can be far away from consistency in selection. For example, with two nested models, if both models are correct, then reasonable selection rules will be consistent in estimation but not necessarily so in selection. See Wegkamp [30] for risk bounds for a modified CV (with an extra complexity penalty) for estimating the regression function.

### 3.1. *Definitions and conditions.*

We first give some useful definitions for comparing estimators in terms of probability.

Let $L(\theta, \widehat{\theta})$ be a loss function. Consider two estimation procedures $\delta$ and $\delta'$ for estimating a parameter $\theta$. Let $\{\widehat{\theta}_{n,1}\}_{n=1}^{\infty}$ and $\{\widehat{\theta}_{n,2}\}_{n=1}^{\infty}$ be the corresponding estimators when applying the two procedures at sample sizes 1, $2, \dots,$ respectively.

Definition 1. Procedure $\delta$ (or $\{\widehat{\theta}_{n,1}\}_{n=1}^{\infty}$, or simply $\widehat{\theta}_{n,1}$) is asymptotically better than $\delta'$ (or $\{\widehat{\theta}_{n,2}\}_{n=1}^{\infty}$, or $\widehat{\theta}_{n,2}$) under the loss function $L(\theta, \widehat{\theta})$ if for every $0 < \epsilon < 1$, there exists a constant $c_\epsilon > 0$ such that when $n$ is large enough,

$$P(L(\theta, \widehat{\theta}_{n,2}) \geq (1 + c_\epsilon) L(\theta, \widehat{\theta}_{n,1})) \geq 1 - \epsilon. \tag{3.1}$$

Remarks. 1. Suppose that $\widehat{\theta}_{n,1}$ and $\widetilde{\theta}_{n,1}$ are asymptotically equivalent under the loss function $L(\theta, \widehat{\theta})$ in the sense that $L(\theta, \widehat{\theta}_{n,1})/L(\theta, \widetilde{\theta}_{n,1}) \to 1$ in probability. If $\widehat{\theta}_{n,1}$ is asymptotically better than $\widehat{\theta}_{n,2}$, then obviously $\widetilde{\theta}_{n,1}$ is also asymptotically better than $\widehat{\theta}_{n,2}$.

2. It seems clear that to evaluate a selection method that chooses between two procedures, the procedures need to be rankable. When two procedures are asymptotically equivalent, one may need to examine finer differences (e.g., higher order behavior) to compare them. When two procedures can have "ties" in the sense that on a set with a nonvanishing probability the two procedures are identical or behave the same, it becomes tricky to define consistency of selection generally speaking.

The requirement in (3.1) is sensible. Obviously, if $L(\theta,\widehat{\theta}_{n,2})/L(\theta,\widehat{\theta}_{n,1}) \to \infty$ in probability, by the definition, $\delta$ is asymptotically better than $\delta'$ under $L$. The concept is also useful for comparing procedures that converge at the same order. In particular, it is worth pointing out that in the context of linear regression with an appropriate loss (e.g., global $L_2$ loss), when two correct models are compared, typically the estimator based on the one with a lower dimension is asymptotically better. As is expected, in the same context, an incorrect subset model yields an asymptotically worse estimator than any one based on a correct model.

We next define consistency of a model/procedure selection rule from the perspective of selecting the best (or better when there are only two candidates) model/procedure. Unless otherwise stated, consistency of CV refers to this notion in the rest of the paper.

Definition 2. Assume that one of the candidate regression procedures, say $\delta^*$, is asymptotically better than the other candidate procedures. A selection rule is said to be consistent if the probability of selecting $\delta^*$ approaches 1 as $n \to \infty$.

This concept of consistency for a selection rule is more general than the definition that a selection rule is consistent if the true model (when existing and being considered) is selected with probability approaching 1. Obviously the latter does not apply for comparing two general regression procedures.

Let $\{a_n\}$ be a sequence of positive numbers approaching zero. The following simple definition concerns the rate of convergence in probability (cf. Stone [23]).

Definition 3. A procedure $\delta$ (or $\{\widehat{\theta}_n\}_{n=1}^{\infty}$) is said to converge exactly at rate $\{a_n\}$ in probability under the loss $L$ if $L(\theta,\widehat{\theta}_n) = O_p(a_n)$, and for every $0 < \epsilon < 1$, there exists $c_\epsilon > 0$ such that when $n$ is large enough, $P(L(\theta,\widehat{\theta}_n) \geq c_\epsilon a_n) \geq 1-\epsilon$.

Clearly, the latter part of the condition in the above definition says that the estimator does not converge faster than the rate $a_n$.

Define the $L_q$ norm

$$\|f\|_q = \begin{cases} \left(\int |f(x)|^q P_X(dx)\right)^{1/q}, & \text{for } 1 \leq q < \infty, \\ \operatorname{ess\,sup} |f|, & \text{for } q = \infty, \end{cases}$$

where $P_X$ denotes the probability distribution of $X_1$.

We assume that $\widehat{f}_{n,1}$ converges exactly at rate $p_n$ in probability and $\widehat{f}_{n,2}$ converges exactly at rate $q_n$ in probability under the $L_2$ loss. Note that $p_n$ and $q_n$ may or may not converge at the same rate.

CONDITION 0 (*Error variances*). The error variances $E(\varepsilon_i^2|X_i)$ are upper bounded by a constant $\overline{\sigma}^2 > 0$ almost surely for all $i \geq 1$.

This is a mild condition on the errors, which does not require them to be identically distributed. We also need to control the sup-norm of the estimators.

CONDITION 1 (*Sup-norm of the estimators*). There exists a sequence of positive numbers $A_n$ such that for $j = 1, 2$, $\|f - \widehat{f}_{n,j}\|_\infty = O_p(A_n)$.

The condition almost always holds. But for our main result to be helpful, the constants $A_n$ need to be suitably controlled.

We mention some useful sufficient conditions that imply Condition 1. One is that for $j = 1, 2$, $\|f - \widehat{f}_{n,j}\|_\infty$ is bounded in probability. Another stronger condition is that $E\|f - \widehat{f}_{n,j}\|_\infty$, is uniformly bounded in $n$ for $j = 1, 2$. It clearly holds if with probability 1, $\|f - \widehat{f}_{n,j}\|_\infty$ is upper bounded by a constant $A > 0$, which is satisfied if the true regression function is bounded between two known constants and the estimators are accordingly restricted.

In order to have consistency in selection, we need that one procedure is better than the other.

CONDITION 2 (*One procedure being better*). Under the $L_2$ loss, either $\delta_1$ is asymptotically better than $\delta_2$, or $\delta_2$ is asymptotically better than $\delta_1$.

Clearly there are situations where neither of two competing procedures is asymptotically better than the other. In such a case, the concept of consistency is hard to define (or may be irrelevant). Note that under Condition 2, if $\delta_1$ is asymptotically better than $\delta_2$, then we have $p_n = O(q_n)$. Clearly, if there exists a constant $C > 1$ such that with probability approaching 1, $\|f - \widehat{f}_{n,2}\|_2 \geq C\|f - \widehat{f}_{n,1}\|_2$ or even $\|f - \widehat{f}_{n,2}\|_2 / \|f - \widehat{f}_{n,1}\|_2 \to \infty$ in probability, then Condition 2 is satisfied with $\delta_1$ being better.

Another quantity, namely, $\|f - \widehat{f}_{n,j}\|_4 / \|f - \widehat{f}_{n,j}\|_2$, is involved in our analysis. Obviously, the ratio is lower bounded by 1 and there cannot be any general positive upper bound. For various estimators, the ratio can be controlled. We need the following condition.

CONDITION 3 (*Relating $L_4$ and $L_2$ losses*). There exists a sequence of positive numbers $\{M_n\}$ such that for $j = 1, 2$, $\|f - \widehat{f}_{n,j}\|_4 / \|f - \widehat{f}_{n,j}\|_2 = O_p(M_n)$.

For many familiar infinite-dimensional classes of regression functions, the optimal rates of convergence under $L_4$ and $L_2$ are the same. If we consider an optimal estimator (in rate) under $L_4$, then we can take $M_n$ to be 1 for a typical $f$ in such a function class. Parametric estimators typically have the ratio upper bounded in probability (i.e., $M_n = 1$) under some mild conditions.

For some nonparametric estimators, the sup-norm risk is often of only a slightly higher order than that under $L_p$ for $p < \infty$. For example, for Hölder classes $\Sigma(\beta, L) = \{f : |f^{(m)}(x) - f^{(m)}(y)| \le L|x - y|^\alpha\}$, where $m = \lfloor \beta \rfloor$ is an integer, $0 < \alpha \le 1$ and $\alpha = \beta - m$; also $\|f\|_\infty$ is bounded, or for Sobolev classes, the rates of convergence under the sup-norm distance and $L_p$ ($p < \infty$) are different only by a logarithmic factor (see, e.g., Stone [24] and Nemirovski [16]). If one takes an optimal or near-optimal estimator under the $L_\infty$ loss, Condition 3 is satisfied typically with $M_n$ being a logarithmic term.

3.2. *The main theorem.* Let $I^* = 1$ if $\delta_1$ is asymptotically better than $\delta_2$ and $I^* = 2$ if $\delta_2$ is asymptotically better than $\delta_1$. Let $\widehat{I}_n = 1$ if $\mathrm{CV}(\widehat{f}_{n_1,1}) \le \mathrm{CV}(\widehat{f}_{n_1,2})$ and otherwise $\widehat{I}_n = 2$.

THEOREM 1. *Under Conditions 0–3, if the data splitting satisfies (1) $n_2 \to \infty$ and $n_1 \to \infty$; (2) $n_2 M_{n_1}^{-4} \to \infty$; and (3)*

$$\sqrt{n_2}\max(p_{n_1}, q_{n_1})/(1 + A_{n_1}) \to \infty,$$

*then the delete-$n_2$ CV is consistent, that is, $P(\widehat{I}_n \ne I^*) \to 0$ as $n \to \infty$.*

REMARKS. 1. The third requirement above is equivalent to $\sqrt{n_2}\max(p_{n_1}, q_{n_1}) \to \infty$ and $\sqrt{n_2}\max(p_{n_1}, q_{n_1})/A_{n_1} \to \infty$. The latter has no effect, for example, when the estimators being compared by CV both converge in the sup-norm in probability. The effect of $M_{n_1}$ is often more restrictive than $A_{n_1}$ and sometimes can be more complicated to deal with.

2. Theorem 1 is stated under a single data generating distribution. Obviously, model selection becomes useful when there are various possible data generating mechanisms (e.g., corresponding to different parametric or nonparametric families of regression functions or different assumptions on the errors) that are potentially suitable for the data at hand. In the linear regression context, a number of finite-dimensional models is considered, and in the literature a model selection rule is said to be consistent when the true model is selected with probability going to 1 no matter what the true model is (assumed to be among the candidates). Clearly our theorem can be used to get such a result when multiple scenarios of the data generating process are possible. To that end, one just needs to find a data splitting ratio that works for all the scenarios being considered.

3. A potentially serious disadvantage of cross validation is that when two candidate regression procedures are hard to distinguish, the forced action of choosing a single winner can substantially damage the accuracy of estimating the regression function. An alternative is to average the estimates. See Yang [33, 34] for references and theoretical results on combining models/procedures and simulation results that compare CV and a model combining the procedure Adaptive Regression by Mixing (ARM).

It should be pointed out that the conclusion is sharp in the sense that there are cases in which the sufficient conditions on data splitting are also necessary. See Section 3.4 for details.

The "ideal" norm conditions are that $M_n = O(1)$ and $A_n = O(1)$. Since almost always $p_n$ and $q_n$ are at least of order $n^{-1/2}$, the most stringent third requirement then is $n_2/n_1 \to \infty$.

COROLLARY 1. *Under the same conditions as in Theorem 1, if $M_n = O(1)$ and $A_n = O(1)$, then the delete-$n_2$ CV is consistent for each of the following two cases:*

(i) $\max(p_n, q_n) = O(n^{-1/2})$, *with the choice $n_1 \to \infty$ and $n_2/n_1 \to \infty$;*
(ii) $\max(p_n, q_n)n^{1/2} \to \infty$, *with any choice such that $n_1 \to \infty$ and $n_1/n_2 = O(1)$.*

From the corollary, if the $L_q$ norms of $\widehat{f}_{n,1} - f$ and $\widehat{f}_{n,2} - f$ behave "nicely" for $q = 2, 4$ and $\infty$, then when at least one of the regression procedures converges at a nonparametric rate under the $L_2$ loss, any split ratio in CV works for consistency as long as both sizes tend to infinity and the estimation size is no bigger (in order) than the evaluation size. Note that this splitting requirement is sufficient but not necessary. For example, if $\max(p_n, q_n) = O(n^{-1/4})$, then the condition $\sqrt{n_2}\max(p_{n_1}, q_{n_1}) \to \infty$ becomes $n_2/n_1^{1/2} \to \infty$ [i.e., $n_1 = o(n_2^2)$]. Thus, for example, we can take $n_1 = n - \lfloor\sqrt{n}\log n\rfloor$ and $n_2 = \lfloor\sqrt{n}\log n\rfloor$, in which case the estimation proportion is dominating. This is in sharp contrast to the requirement of a much larger evaluation size ($n_2/n_1 \to \infty$) when both of the regression procedures converge at the parametric rate $n^{-1}$, which was discovered by Shao [19] in the context of linear regression with fixed design. From Shao's results, one may expect that when $p_n$ and $q_n$ are of the same order, the condition $n_2/n_1 \to \infty$ may also be needed for consistency more generally. But, interestingly, our result shows that this is not the case. Thus there is a paradigm shift in terms of splitting proportion for CV when at least one of the regression procedures is nonparametric.

The result also suggests that in general delete-1 (or delete a small fraction) is not geared toward finding out which candidate procedure is the better one.

Note that for various estimators based on local averaging or series expansion, we do not necessarily have the "ideal" norm requirement in Condition 3 met with $M_n = O(1)$, but may have $\|f - \widehat{f}_n\|_4 \leq \|f - \widehat{f}_n\|_\infty \leq a_n \|f - \widehat{f}_n\|_2$ for some deterministic sequence $a_n$ (possibly converging to $\infty$ at a polynomial order $n^\gamma$ with $\gamma > 0$). Then for applying the theorem, we may need $n_2/n_1^{4\gamma} \to \infty$. This may or may not add any further restriction on the data splitting ratio in CV beyond the other requirements in the theorem, depending on the value of $\gamma$ in relation to $p_n, q_n$ and $A_n$.

3.3. *CV with multiple data splittings.* For Theorem 1, the data splitting in CV is done only once (and hence the name *cross* validation may not be appropriate there). Clearly, the resulting estimator depends on the order of the observations. In real applications, one may do any of the following: (1) consider all possible splits with the same ratio (this is called multifold CV; see, e.g., Zhang [35]); (2) the same as in (1) but consider only a sample of all possible splits (this is called repeated learning-testing; see, e.g., Burman [4]); (3) divide the data into $r$ subgroups and do prediction one at a time for each subgroup based on estimation using the rest of the subgroups (this is called $r$-fold CV; see Breiman, Friedman, Olshen and Stone [3]). When multiple splittings are used, there are two natural ways to proceed. One is to first average the prediction errors over the different splittings and then select the procedure that minimizes the average prediction error. Another is to count the number of times each candidate is preferred under the different splittings and then the candidate with the highest count is the overall winner. Such a voting is natural to consider for model selection. To make a distinction, we call the former (i.e., CV with averaging) CV-a and the latter (i.e., CV with voting) CV-v. When focusing on linear models with fixed design, theoretical properties for multifold CV-a or $r$-fold CV-a were derived under assumptions on the design matrix (see, e.g., Zhang [35] and Shao [19]). We next show that under the same conditions used for a single data splitting in the previous subsection, CV-v based on multiple data splittings is also consistent in selection when the observations are identically distributed.

Let $\pi$ denote a permutation of the observations. Let $\mathrm{CV}_\pi(\widehat{f}_{n_1,j})$ be the criterion value as defined in (2.1) except that the data splitting is done after the permutation $\pi$. If $\mathrm{CV}_\pi(\widehat{f}_{n_1,1}) \leq \mathrm{CV}_\pi(\widehat{f}_{n_1,2})$, then let $\tau_\pi = 1$ and otherwise let $\tau_\pi = 0$.

Let $\Pi$ denote the set of all $n!$ permutations of the observations. If $\sum_{\pi \in \Pi} \tau_\pi \geq \frac{n!}{2}$, then we select $\delta_1$, and otherwise $\delta_2$ is selected.

THEOREM 2. *Under the conditions in Theorem 1 and that the observations are independent and identically distributed, the CV-v procedure above is consistent in selection.*

PROOF. Without loss of generality, assume that $\delta_1$ is better than $\delta_2$. Let $W$ denote the values of $(X_1, Y_1), \dots, (X_n, Y_n)$ (ignoring the orders). Under the i.i.d. assumption on the observations, obviously, conditional on $W$, every ordering of these values has exactly the same probability and thus $P(\mathrm{CV}(\widehat{f}_{n_1,1}) \le \mathrm{CV}(\widehat{f}_{n_1,2}))$ is equal to

$$EP(\mathrm{CV}(\widehat{f}_{n_1,1}) \le \mathrm{CV}(\widehat{f}_{n_1,2})|W) = E\left(\frac{\sum_{\pi\in\Pi}\tau_\pi}{n!}\right).$$

From Theorem 1, under the given conditions, $P(\mathrm{CV}(\widehat{f}_{n_1,1}) \le \mathrm{CV}(\widehat{f}_{n_1,2})) \to 1$. Thus $E(\sum_{\pi\in\Pi}\tau_\pi/n!) \to 1$. Since $\sum_{\pi\in\Pi}\tau_\pi/n!$ is between 0 and 1, for its expectation to converge to 1, we must have $\sum_{\pi\in\Pi}\tau_\pi/n! \to 1$ in probability. Consequently $P(\sum_{\pi\in\Pi}\tau_\pi \ge n!/2) \to 1$. This completes the proof of Theorem 2. $\square$

From the above proof, it is clearly seen that the consistency result also holds for the aforementioned voting-based repeated learning-testing and $r$-fold CV-v methods.

Note that for the CV-a methods, for each candidate procedure, we average $\mathrm{CV}(\widehat{f}_{n_1,j})$ over the different data splittings first and then compare the criterion values to select a winner. Intuitively, since the CV-v only keeps the ranking of the procedures for each splitting, it may have lost some useful information in the data. If so, the CV-v may perform worse than CV-a. However, in terms of the consistency property in selection, it is unlikely that the two versions of CV are essentially different. Indeed, if $P(\sum_{\pi\in\Pi}(\mathrm{CV}_\pi(\widehat{f}_{n_1,2}) - \mathrm{CV}_\pi(\widehat{f}_{n_1,1})) > 0) \to 1$, due to symmetry, one expects that for the majority of splittings in $\Pi$, with high probability, we have $\mathrm{CV}_\pi(\widehat{f}_{n_1,2}) - \mathrm{CV}_\pi(\widehat{f}_{n_1,1}) > 0$. Then $P(\sum_{\pi\in\Pi} I_{(\mathrm{CV}_\pi(\widehat{f}_{n_1,2}) - \mathrm{CV}_\pi(\widehat{f}_{n_1,1})>0)} > \frac{|\Pi|}{2})$ is close to 1.

Based on the above reasoning, we conjecture that the two CV methods generally share the same status of consistency in selection (i.e., if one is consistent so is the other). This of course does not mean that they perform similarly at a finite sample size. In the simulations reported in Section 4, we see that for comparing two converging parametric procedures, CV-v performs much worse than CV-a, but when CV-a selects the best procedure with probability closer to 1, their difference becomes small; for comparing a converging nonparametric procedure with a converging parametric one, CV-v often performs better. We tend to believe that the differences between CV-v and CV-a are second-order effects, which are hard to quantify in general.

3.4. *Are the sufficient conditions on data splitting also necessary for consistency*? Theorem 2 shows that if a single data splitting ensures consistency, voting based on multiple splittings also works. In the reverse direction, one may wonder if multiple splittings with averaging can rescue an inconsistent CV selection with only one splitting. We give a counterexample below.

In the following example, we consider two simple models which allow us to exactly identify the sufficient and necessary conditions for ensuring consistency in selection by a CV-a method. Model 1 is $Y_i = \varepsilon_i$, $i = 1, \ldots, n$, where $\varepsilon_i$ are i.i.d. normal with mean zero and variance $\sigma^2$. Model 2 is $Y_i = \mu + \varepsilon_i$, $i = 1, \ldots, n$, with the same conditions on the errors. Clearly, model 1 is a submodel of model 2. Under model 1, obviously $\widehat{f}_{n,1}(x) = 0$, and under model 2, the maximum likelihood method yields $\widehat{f}_{n,2}(x) = \overline{Y}$. We consider the multifold CV with splitting ratio $n_1 : n_2$, where again $n_1$ is the estimation sample size and $n_2$ is the evaluation size. Note that for the multifold CV-a method, we consider all possible data splittings at the given ratio and the averaging of the prediction errors is done over the splittings. Let $S$ denote the observations in the estimation set. Then

$$\mathrm{CV}(j) = \sum_S \left( \sum_{i \in S^c} (Y_i - \widehat{Y}_{i,j})^2 \right),$$

where $S$ is over all possible data splittings with the given ratio, $S^c$ denotes the complement of $S$, $\widehat{Y}_{i,1} = 0$ and $\widehat{Y}_{i,2} = \frac{1}{n_1} \sum_{l \in S} Y_l$.

Proposition 1. *A sufficient and necessary condition for the above CV to be consistent in selection is that the data splitting ratio satisfies* (1) $n_1 \to \infty$, (2) $n_2/n_1 \to \infty$.

Note that the conditions in this proposition match those sufficient conditions in Theorem 1. Consequently, we know that for the consistency property, multiple splittings (or even all possible splittings) do not help in this case. Although it seems technically difficult to derive a similar general result, we tend to believe that the example represents the typical situation of comparing nested parametric models.

In the context of variable selection in linear regression, Zhang [35] showed that for any fixed splitting ratio, none of multifold CV, repeated learning-testing and $r$-fold CV based on averaging is consistent in selection.

Proof of Proposition 1. With some calculations, under the larger model, we have $\mathrm{CV}(1) = \binom{n-1}{n_2-1}(n\mu^2 + 2\mu \sum_{i=1}^n \varepsilon_i + \sum_{i=1}^n \varepsilon_i^2)$ and CV(2) equals

$$\binom{n-1}{n_2-1} \sum_{i=1}^n \varepsilon_i^2 - 2 \binom{n-2}{n_1-1} \frac{1}{n_1} \sum_{i=1}^n \varepsilon_i \left( \sum_{j \neq i}^n \varepsilon_j \right) + \frac{n_2}{n_1^2} \left( \binom{n-1}{n_1-1} \sum_{i=1}^n \varepsilon_i^2 + 2 \binom{n-2}{n_1-2} \sum_{1 \le j_1 < j_2 \le n}^n \varepsilon_{j_1} \varepsilon_{j_2} \right).$$

Then with more simplifications, we get that CV(1) − CV(2) equals

$$\binom{n-1}{n_2-1}\left(n\mu^2+2\mu\sum_{i=1}^{n}\varepsilon_i\right)$$
$$+\sum_{i=1}^{n}\varepsilon_i^2\cdot\frac{(n-2)!(n_1+1)}{n_1 n_1!(n_2-1)!}-\left(\sum_{i=1}^{n}\varepsilon_i\right)^2\cdot\frac{(n-2)!(n_1+n)}{n_1 n_1!(n_2-1)!}.$$

When $\mu\neq 0$, it is not hard to show that $\mathrm{CV}(1)-\mathrm{CV}(2)$ is positive with probability tending to 1 if and only if $n_1\to\infty$. When model 1 holds, that is, $\mu=0$, $\mathrm{CV}(1)-\mathrm{CV}(2)>0$ is equivalent to $\frac{(n-1)n_1}{n}(\sum_{i=1}^{n}\varepsilon_i)^2>(n+n_1)\sum_{i=1}^{n}(\varepsilon_j-\overline{\varepsilon})^2$. Since $\sum_{i=1}^{n}\varepsilon_i$ is independent of $\sum_{i=1}^{n}(\varepsilon_j-\overline{\varepsilon})^2$ and they have normal and chi-square distributions, respectively, we know that for the probability of the event to go to zero, we must have $n_1/n\to 0$, which is also sufficient. This completes the proof of the proposition. □

**4. Simulation.** In this section, we present some simulation results that are helpful to understand the differences of several versions of CV and the effect of the data splitting ratio.

We consider estimating a regression function on $[0,1]$ with three competing regression methods. The true model is $Y_i=f(X_i)+\varepsilon_i$, $1\leq i\leq n$, where $X_i$ are i.i.d. uniform in the unit interval, and the errors are independent of $X_i$ and are i.i.d. normal with mean zero and standard deviation $\sigma=0.3$. The true function is taken to be one of the three functions

(Case 1) $\quad f_1(x)=1+x,$

(Case 2) $\quad f_2(x)=1+x+0.7(x-0.5)^2,$

(Case 3) $\quad f_3(x)=1+x-\exp(-200(x-0.25)^2).$

In all these cases, for $\sigma$ not too large, a linear trend in the scatter plot is more or less obvious, but when $\sigma$ is not small, it is usually not completely clear whether the true function is simply linear or not. We consider three regression methods: simple linear regression, quadratic regression and smoothing spline. The simulation was conducted using R, where a smoothing spline method is provided. We take the default choice of GCV for smoothing parameter selection.

Clearly, in Case 1, the linear regression is the winner; the quadratic regression is the right one in Case 2 (when the sample size is reasonably large); and the smoothing spline method is the winner in Case 3. Case 1 is the most difficult in comparing the methods because all of the three estimators converge to the true regression function with two converging at the same parametric rate. Case 2 is easier, where the simple linear estimator no longer converges and consequently the task is basically to compare the parametric and the nonparametric estimators (note that the spline estimator converges at a slower order). Case 3 is the easiest, where only the smoothing spline estimator is converging.

We consider the following versions of CV: (1) single splitting CV: the data splitting is done just once; (2) repeated learning-testing (a version of CV-a), denoted as RLT: we randomly split the data into the estimation and evaluation parts 100 times and average the prediction errors over the splittings for each estimator; (3) repeated splittings with voting (a version of CV-v), denoted as RSV: differently from the previous one, we select the best method based on each data splitting and then vote to decide the overall winner. The sample sizes considered are $100, 200, 400, 800$ and 1600. In Case 2, the first one or two sample sizes are not considered for the splitting ratios $3:7$ and $1:9$ because the smoothing spline method has difficulty in parameter estimation due to the small sample size in the estimation part. In Case 3, only the splitting ratios $9:1$ and $5:5$ are included because lower splitting ratios make the CV methods perform perfectly. Note also that the first two sample sizes are not included for the splitting ratio $5:5$ for the same reason as mentioned above.

The results based on 200 replications are presented in Figures 1–3 for the three cases.

From the graphs, we observe the following:

1. In Case 1, at a given splitting ratio, the increase of sample size does not lead to improvement on correct identification of the best estimator. This nicely matches what is expected from Theorem 1: since the simple linear and the quadratic regression estimators both converge at the parametric rate, no matter how large the sample size is, any fixed proportion is not sufficient for consistency in selection.

2. In Case 2, at a given splitting ratio, when the sample size is increased, the ability of CV to identify the best estimator tends to be enhanced. This is also consistent with Theorem 1: in this case, the two converging estimators converge at different rates and thus a fixed splitting ratio is sufficient to ensure consistency in selection.

3. In Case 3, even at splitting ratio $9:1$, RLT and RSV have little difficulty in finding the best estimator. From Theorem 1, in this case, since the smoothing spline estimator is the only converging one, the splitting ratio $n_1:n_2$ is even allowed to go to $\infty$ without sacrificing the property of consistency in selection.

4. The different versions of CV behave quite differently. Overall, RLT seems to be the best, while the single splitting is clearly inferior. In Case 1, for the splitting ratios that favor the estimation size, RSV did poorly. How-ever, with higher and higher splitting ratio toward the evaluation part, it improved dramatically. This is consistent with the understanding that for a sequence of consistent splitting ratios, CV with voting or CV with averaging do not differ asymptotically (although their second-order behavior may be different). In Cases 2 and 3, the CV-v actually performed similarly to or better than the CV-a [the difference is large in Case 2 when the estimation sample size is small (40, 50, 60)].

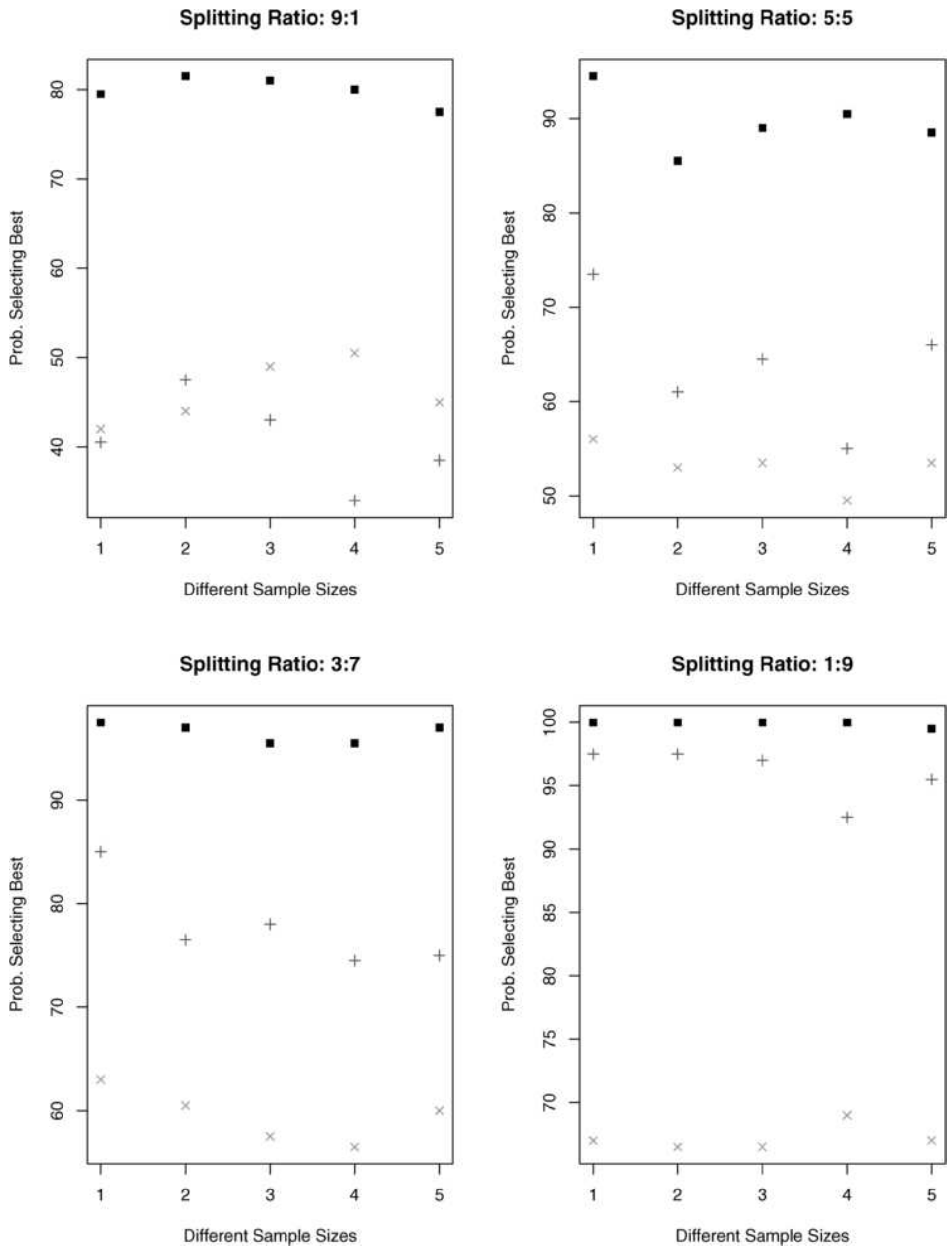

FIG. 1. *Probability of selecting the best estimator for* Case 1. *Coding:* ■ = *repeated learning-testing,* + = *repeated sampling with voting,* × = *single splitting. The sample sizes are* 100, 200, 400, 800 *and* 1600.

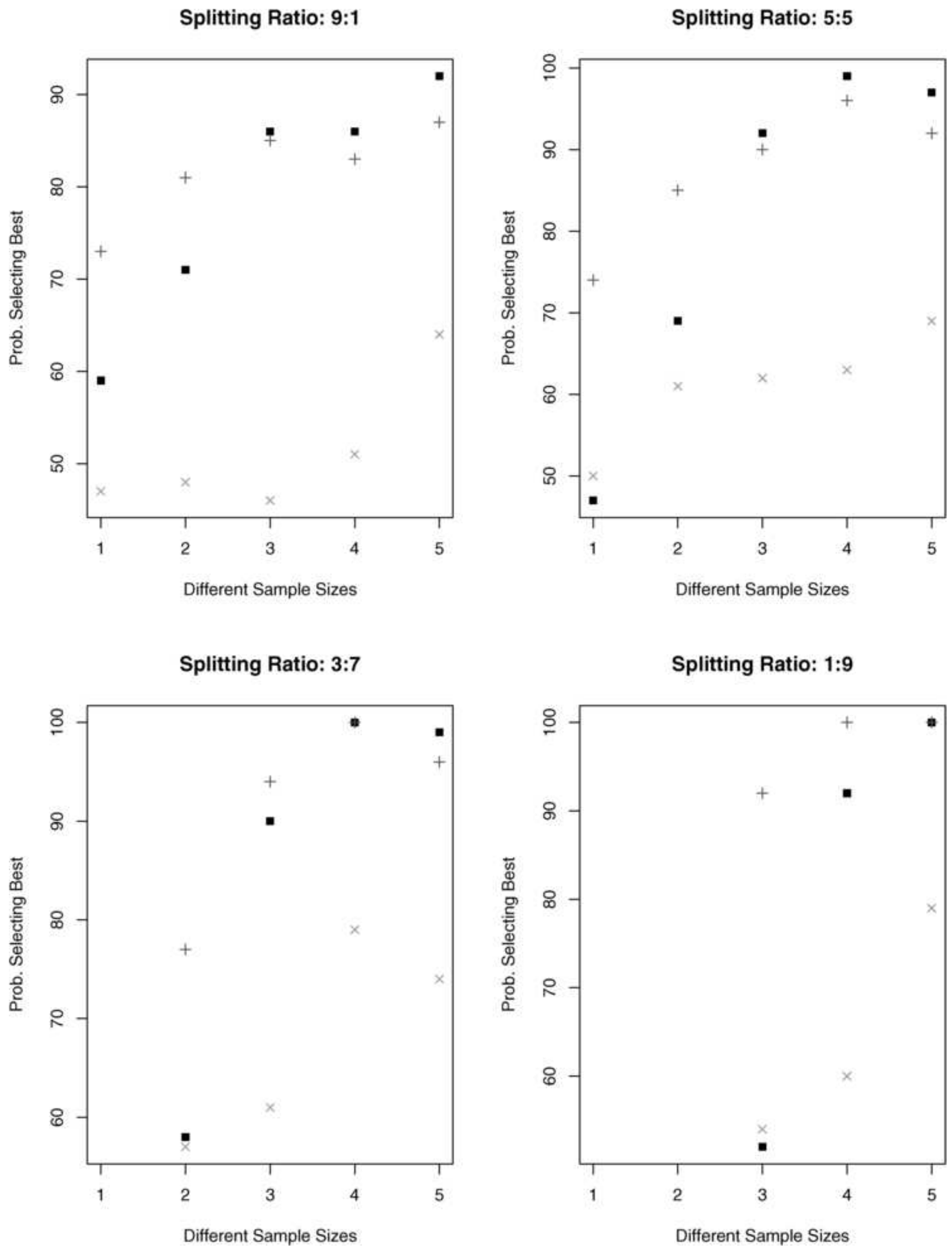

Fig. 2. *Probability of selecting the best estimator for* Case 2. *Coding:* ■ = *repeated learning-testing,* + = *repeated sampling with voting,* × = *single splitting. The sample sizes are* 100, 200, 400, 800 *and* 1600.

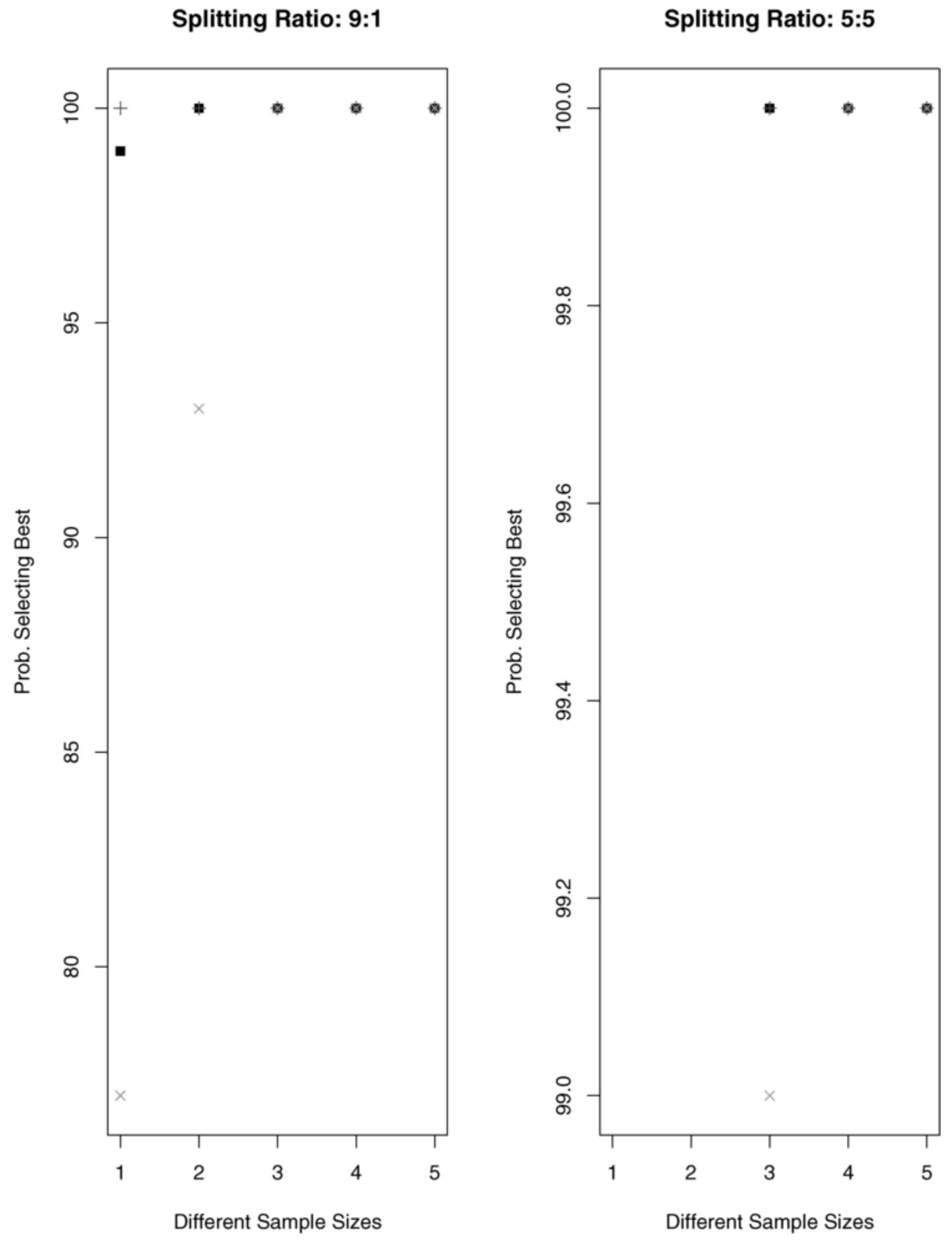

FIG. 3. *Probability of selecting the best estimator for* Case 3. *Coding:* ■ = *repeated learning-testing,* + = *repeated sampling with voting,* × = *single splitting. The sample sizes are* 100, 200, 400, 800 *and* 1600.

In summary, the simulation results are very much in line with the understanding in Section 3.

**5. Concluding remarks.** We have shown that under some sensible conditions on the $L_2$, $L_4$ and $L_\infty$ norms of $\widehat{f} - f$ for the competing estimators, with an appropriate splitting ratio of the data for cross validation, the better model/procedure will be selected with probability converging to 1. Unlike the previous results on CV that focus either on comparing only parametric regression models or on selecting a smoothing parameter in nonparametric regression, our result can be applied generally to compare both parametric and nonparametric models/procedures.

The result also reveals some interesting behavior of CV. Differently from the parametric model selection case, for comparing two models converging at the same nonparametric rate, it is not necessary for the evaluation size in data splitting to be dominating. Actually, the proportion of the evaluation part can even be of a smaller order than the estimation proportion without damaging the property of consistency in selection. An implication is that it may be desirable to take the characteristics of the regression estimators into consideration for data splitting, which to our knowledge has not been seriously addressed in the literature. Based on our result, for comparing two estimators with at least one being nonparametric, half–half splitting is a good choice.

Delete-1 CV is a popular choice in applications. This is usually suitable for estimating the regression function. However, when one's goal is to find which method is more accurate for the data at hand, the proportion of evaluation needs to be much larger. Further research on practical methods for properly choosing the data splitting proportion can be valuable for successful applications of cross validation.

**6. Proof of Theorem 1.** Without loss of generality, assume that $\widehat{f}_{n,1}$ is the asymptotically better estimator by Condition 2. Note that

$$\begin{aligned}\mathrm{CV}(\widehat{f}_{n_1,j}) &= \sum_{i=n_1+1}^{n} (f(X_i) - \widehat{f}_{n_1,j}(X_i) + \varepsilon_i)^2 \\ &= \sum_{i=n_1+1}^{n} \varepsilon_i^2 + \sum_{i=n_1+1}^{n} (f(X_i) - \widehat{f}_{n_1,j}(X_i))^2\end{aligned}$$

$$+ 2 \sum_{i=n_1+1}^{n} \varepsilon_i (f(X_i) - \widehat{f}_{n_1,j}(X_i)).$$

Define for $j = 1, 2$,

$$L_j = \sum_{i=n_1+1}^{n} (f(X_i) - \widehat{f}_{n_1,j}(X_i))^2 + 2 \sum_{i=n_1+1}^{n} \varepsilon_i (f(X_i) - \widehat{f}_{n_1,j}(X_i)).$$

Then $\mathrm{CV}(\widehat{f}_{n_1,1}) \le \mathrm{CV}(\widehat{f}_{n_1,2})$ is equivalent to $L_1 \le L_2$ and thus also equivalent to

$$\begin{aligned} 2 \sum_{i=n_1+1}^{n} \varepsilon_i (\widehat{f}_{n_1,2}(X_i) - \widehat{f}_{n_1,1}(X_i)) \le & \sum_{i=n_1+1}^{n} (f(X_i) - \widehat{f}_{n_1,2}(X_i))^2 \\ & - \sum_{i=n_1+1}^{n} (f(X_i) - \widehat{f}_{n_1,1}(X_i))^2. \end{aligned}$$

Conditional on $Z^1$ and $X^2 = (X_{n_1+1}, \dots, X_n)$, assuming $\sum_{i=n_1+1}^{n} (f(X_i) - \widehat{f}_{n_1,2}(X_i))^2$ is larger than $\sum_{i=n_1+1}^{n} (f(X_i) - \widehat{f}_{n_1,1}(X_i))^2$, by Chebyshev's inequality, we have

$$\begin{aligned} & P(\mathrm{CV}(\widehat{f}_{n_1,1}) > \mathrm{CV}(\widehat{f}_{n_1,2}) | Z^1, X^2) \\ & \quad \le \min\Biggl(1, 4\overline{\sigma}^2 \sum_{i=n_1+1}^{n} (\widehat{f}_{n_1,2}(X_i) - \widehat{f}_{n_1,1}(X_i))^2 \\ & \qquad \times \Biggl[\Biggl( \sum_{i=n_1+1}^{n} (f(X_i) - \widehat{f}_{n_1,2}(X_i))^2 \\ & \qquad\quad - \sum_{i=n_1+1}^{n} (f(X_i) - \widehat{f}_{n_1,1}(X_i))^2 \Biggr)^2 \Biggr]^{-1} \Biggr). \end{aligned}$$

Let $Q_n$ denote the ratio in the upper bound in the above inequality and let $S_n$ be the event of $\sum_{i=n_1+1}^{n} (f(X_i) - \widehat{f}_{n_1,2}(X_i))^2 > \sum_{i=n_1+1}^{n} (f(X_i) - \widehat{f}_{n_1,1}(X_i))^2$. It follows that

$$\begin{aligned} & P(\mathrm{CV}(\widehat{f}_{n_1,1}) > \mathrm{CV}(\widehat{f}_{n_1,2})) \\ & \quad = P(\{\mathrm{CV}(\widehat{f}_{n_1,1}) > \mathrm{CV}(\widehat{f}_{n_1,2})\} \cap S_n) + P(\{\mathrm{CV}(\widehat{f}_{n_1,1}) > \mathrm{CV}(\widehat{f}_{n_1,2})\} \cap S_n^c) \\ & \quad \le E(P(\mathrm{CV}(\widehat{f}_{n_1,1}) > \mathrm{CV}(\widehat{f}_{n_1,2}) | Z^1, X^2) I_{S_n}) + P(S_n^c) \\ & \quad \le E \min(1, Q_n) + P(S_n^c). \end{aligned}$$

If we can show that $P(S_n^c) \to 0$ and $Q_n \to 0$ in probability as $n \to \infty$, then due to the boundedness of $\min(1, Q_n)$ [which implies that the random variables $\min(1, Q_n)$ are uniformly integrable], we have $P(\mathrm{CV}(\widehat{f}_{n_1,1}) >$ $\mathrm{CV}(\widehat{f}_{n_1,2})$) converges to zero as $n \to \infty$. Suppose we can show that for every $\epsilon > 0$, there exists $\alpha_\epsilon > 0$ such that when $n$ is large enough,

$$P\left(\frac{\sum_{i=n_1+1}^{n}(f(X_i)-\widehat{f}_{n_1,2}(X_i))^2}{\sum_{i=n_1+1}^{n}(f(X_i)-\widehat{f}_{n_1,1}(X_i))^2} \geq 1+\alpha_\epsilon\right) \geq 1-\epsilon. \tag{6.1}$$

Then $P(S_n) \geq 1-\epsilon$ and thus $P(S_n^c) \to 0$ as $n \to \infty$. By the triangle inequality,

$$\begin{aligned}&\sum_{i=n_1+1}^{n}(\widehat{f}_{n_1,2}(X_i)-\widehat{f}_{n_1,1}(X_i))^2\\ &\quad\leq 2\sum_{i=n_1+1}^{n}(f(X_i)-\widehat{f}_{n_1,1}(X_i))^2 + 2\sum_{i=n_1+1}^{n}(f(X_i)-\widehat{f}_{n_1,2}(X_i))^2.\end{aligned}$$

Then with probability no less than $1-\epsilon$, $Q_n$ is upper bounded by

$$\begin{aligned}&\frac{8\overline{\sigma}^2(\sum_{i=n_1+1}^{n}(f(X_i)-\widehat{f}_{n_1,1}(X_i))^2+\sum_{i=n_1+1}^{n}(f(X_i)-\widehat{f}_{n_1,2}(X_i))^2)}{((1-1/(1+\alpha_\epsilon))\sum_{i=n_1+1}^{n}(f(X_i)-\widehat{f}_{n_1,2}(X_i))^2)^2}\\ &\quad\leq \frac{8\overline{\sigma}^2(1+1/(1+\alpha_\epsilon))}{(1-1/(1+\alpha_\epsilon))^2\sum_{i=n_1+1}^{n}(f(X_i)-\widehat{f}_{n_1,2}(X_i))^2}.\end{aligned} \tag{6.2}$$

From (6.1) and (6.2), to show $P(S_n^c) \to 0$ and $Q_n \to 0$ in probability, it suffices to show (6.1) and

$$\sum_{i=n_1+1}^{n}(f(X_i)-\widehat{f}_{n_1,2}(X_i))^2 \to \infty \qquad \text{in probability.} \tag{6.3}$$

Suppose a slight relaxation of Condition 1 holds: for every $\epsilon > 0$, there exists $A_{n_1,\epsilon}$ such that when $n_1$ is large enough, $P(\|f-\widehat{f}_{n_1,j}\|_\infty \geq A_{n_1,\epsilon}) \leq \epsilon$ for $j = 1,2$. Let $H_{n_1}$ be the event $\{\max(\|f-\widehat{f}_{n_1,1}\|_\infty, \|f-\widehat{f}_{n_1,2}\|_\infty) \leq A_{n_1,\epsilon}\}$. Then on $H_{n_1}$, we have $W_i = (f(X_i)-\widehat{f}_{n_1,j}(X_i))^2 - \|f-\widehat{f}_{n_1,j}\|_2^2$ is bounded between $-(A_{n_1,\epsilon})^2$ and $(A_{n_1,\epsilon})^2$. Notice that conditional on $Z^1$ and $H_{n_1}$,

$$\mathrm{Var}_{Z^1}(W_{n_1+1}) \leq E_{Z^1}(f(X_{n_1+1})-\widehat{f}_{n_1,j}(X_{n_1+1}))^4 = \|f-\widehat{f}_{n_1,j}\|_4^4,$$

where the subscript $Z^1$ in $\mathrm{Var}_{Z^1}$ and $E_{Z^1}$ is used to denote the conditional expectation given $Z^1$. Thus conditional on $Z^1$, on $H_{n_1}$, by Bernstein's inequality (see, e.g., Pollard [18], page 193), for each $x > 0$, we have

$$\begin{aligned}&P_{Z^1}\left(\sum_{i=n_1+1}^{n}(f(X_i)-\widehat{f}_{n_1,1}(X_i))^2 - n_2\|f-\widehat{f}_{n_1,1}\|_2^2 \geq x\right)\\ &\quad\leq \exp\left(-\frac{1}{2}\frac{x^2}{n_2\|f-\widehat{f}_{n_1,1}\|_4^4+(2(A_{n_1,\epsilon})^2x/3)}\right).\end{aligned}$$

Taking $x = \beta_n n_2 \|f - \widehat{f}_{n_1,1}\|_2^2$, the above inequality becomes

$$P_{Z^1}\left(\sum_{i=n_1+1}^{n} (f(X_i) - \widehat{f}_{n_1,1}(X_i))^2 \geq (1+\beta_n) n_2 \|f - \widehat{f}_{n_1,1}\|_2^2\right)$$
$$\leq \exp\left(-\frac{1}{2}\frac{\beta_n^2 n_2 \|f - \widehat{f}_{n_1,1}\|_2^4}{\|f - \widehat{f}_{n_1,1}\|_4^4 + (2(A_{n_1,\epsilon})^2 \beta_n/3)\|f - \widehat{f}_{n_1,1}\|_2^2}\right).$$

Under Condition 2, for every $\epsilon > 0$, there exists $\alpha'_\epsilon > 0$ such that when $n$ is large enough, $P(\|f - \widehat{f}_{n_1,2}\|_2^2 / \|f - \widehat{f}_{n_1,1}\|_2^2 \leq 1 + \alpha'_\epsilon) \leq \epsilon$. Take $\beta_n$ such that $1 + \beta_n = \|f - \widehat{f}_{n_1,2}\|_2^2 / ((1 + \alpha'_\epsilon/2)\|f - \widehat{f}_{n_1,1}\|_2^2)$. Then with probability at least $1 - \epsilon$, $\beta_n \geq (\alpha'_\epsilon/2)/(1 + \alpha'_\epsilon/2)$. Let $D_n$ denote this event and let $S1 = \sum_{i=n_1+1}^{n} (f(X_i) - \widehat{f}_{n_1,1}(X_i))^2$. Then on $D_n$ we have

$$\beta_n \geq \alpha'_\epsilon \|f - \widehat{f}_{n_1,2}\|_2^2 / (2(1+\alpha'_\epsilon)(1+\alpha'_\epsilon/2)\|f - \widehat{f}_{n_1,1}\|_2^2),$$

$$\begin{aligned}
&P_{Z^1}(S1 \geq (1+\beta_n) n_2 \|f - \widehat{f}_{n_1,1}\|_2^2)\\
&\quad = P_{Z^1}\left(S1 \geq \frac{n_2}{1+\alpha'_\epsilon/2}\|f - \widehat{f}_{n_1,2}\|_2^2\right)\\
&\quad \leq P_{Z^1}\left(S1 \geq \left(1 + \frac{\alpha'_\epsilon \|f - \widehat{f}_{n_1,2}\|_2^2}{2(1+\alpha'_\epsilon)(1+\alpha'_\epsilon/2)\|f - \widehat{f}_{n_1,1}\|_2^2}\right) n_2 \|f - \widehat{f}_{n_1,1}\|_2^2\right)\\
&\quad \leq \exp\left(-\frac{(\alpha'_\epsilon)^2}{8(1+\alpha'_\epsilon)^2(1+\alpha'_\epsilon/2)^2}\right.\\
&\qquad\qquad \left.\times \frac{n_2 \|f - \widehat{f}_{n_1,2}\|_2^4}{\|f - \widehat{f}_{n_1,1}\|_4^4 + (\alpha'_\epsilon (A_{n_1,\epsilon})^2/(3(1+\alpha'_\epsilon)(1+\alpha'_\epsilon/2)))\|f - \widehat{f}_{n_1,2})\|_2^2}\right).
\end{aligned}$$

If we have

$$\frac{n_2 \|f - \widehat{f}_{n_1,2}\|_2^4}{\|f - \widehat{f}_{n_1,1}\|_4^4} \to \infty \qquad \text{in probability,} \tag{6.4}$$

$$\frac{n_2 \|f - \widehat{f}_{n_1,2}\|_2^2}{(A_{n_1,\epsilon})^2} \to \infty \qquad \text{in probability,} \tag{6.5}$$

then the upper bound in the last inequality above converges to zero in probability. From these pieces, we can conclude that

$$\begin{aligned}
&P\left(\sum_{i=n_1+1}^{n} (f(X_i) - \widehat{f}_{n_1,1}(X_i))^2 \geq \frac{n_2}{1+\alpha'_\epsilon/2}\|f - \widehat{f}_{n_1,2}\|_2^2\right)\\
&\qquad \leq 3\epsilon + \Delta(\epsilon, n),
\end{aligned} \tag{6.6}$$

for some $\Delta(\epsilon,n)\to 0$ as $n\to\infty$. Indeed, for every given $\epsilon>0$, when $n$ is large enough,

$$\begin{aligned}
&P\biggl(\frac{1}{n_2}\sum_{i=n_1+1}^{n}(f(X_i)-\widehat{f}_{n_1,1}(X_i))^2\ge\frac{1}{1+\alpha'_\epsilon/2}\|f-\widehat{f}_{n_1,2}\|_2^2\biggr)\\
&\quad\le P(H_{n_1}^c)+P(D_n^c)\\
&\qquad+P\biggl(H_{n_1}\cap D_n\cap\biggl\{\frac{1}{n_2}\sum_{i=n_1+1}^{n}(f(X_i)-\widehat{f}_{n_1,1}(X_i))^2\\
&\qquad\qquad\ge\frac{1}{1+\alpha'_\epsilon/2}\|f-\widehat{f}_{n_1,2}\|_2^2\biggr\}\biggr)\\
&\quad\le 3\epsilon+EP\biggl(H_{n_1}\cap D_n\cap\biggl\{\frac{1}{n_2}\sum_{i=n_1+1}^{n}(f(X_i)-\widehat{f}_{n_1,1}(X_i))^2\\
&\qquad\qquad\ge\frac{1}{1+\alpha'_\epsilon/2}\|f-\widehat{f}_{n_1,2}\|_2^2\biggr\}\Big|Z^1\biggr)\\
&\quad\le 3\epsilon+E\exp\biggl(-\frac{(\alpha'_\epsilon)^2}{8(1+\alpha'_\epsilon)^2(1+\alpha'_\epsilon/2)^2}\\
&\qquad\times\frac{n_2\|f-\widehat{f}_{n_1,2}\|_2^4}{\|f-\widehat{f}_{n_1,1}\|_4^4+(\alpha'_\epsilon(A_{n_1,\epsilon})^2/(3(1+\alpha'_\epsilon)(1+\alpha'_\epsilon/2)))\|f-\widehat{f}_{n_1,2}\|_2^2}\biggr)\\
&\quad\triangleq 3\epsilon+\Delta(\epsilon,n),
\end{aligned}$$

where the expectation in the upper bound of the last inequality above [i.e., $\Delta(\epsilon,n)$] converges to zero due to the convergence in probability to zero of the random variables of the exponential expression and their uniform integrability (since they are bounded above by 1), provided that (6.4) and (6.5) hold. The assertion of (6.6) then follows.

For the other estimator, similarly, for $0<\widetilde{\beta}_n<1$, we have

$$\begin{aligned}
&P_{Z^1}\biggl(\sum_{i=n_1+1}^{n}(f(X_i)-\widehat{f}_{n_1,2}(X_i))^2\le(1-\widetilde{\beta}_n)n_2\|f-\widehat{f}_{n_1,2}\|_2^2\biggr)\\
&\qquad\le\exp\biggl(-\frac{1}{2}\frac{n_2\widetilde{\beta}_n^2\|f-\widehat{f}_{n_1,2}\|_2^4}{\|f-\widehat{f}_{n_1,2}\|_4^4+(2(A_{n_1,\epsilon})^2\widetilde{\beta}_n/3)\|f-\widehat{f}_{n_1,2}\|_2^2}\biggr).
\end{aligned}$$

If we have

$$\frac{n_2\widetilde{\beta}_n^2\|f-\widehat{f}_{n_1,2}\|_2^4}{\|f-\widehat{f}_{n_1,2}\|_4^4}\to\infty\qquad\text{in probability,}\tag{6.7}$$

$$\frac{n_2\widetilde{\beta}_n\|f-\widehat{f}_{n_1,2}\|_2^2}{(A_{n_1,\epsilon})^2}\to\infty \qquad \text{in probability,} \tag{6.8}$$

then following a similar argument used for $\widehat{f}_{n_1,1}$, we have

$$P\Bigg(\sum_{i=n_1+1}^{n}(f(X_i)-\widehat{f}_{n_1,2}(X_i))^2\le(1-\widetilde{\beta}_n)n_2\|f-\widehat{f}_{n_1,2}\|_2^2\Bigg)\to 0. \tag{6.9}$$

From this, if $n_2\|f-\widehat{f}_{n_1,2}\|_2^2\to\infty$ in probability and $\widetilde{\beta}_n$ is bounded away from 1, then (6.3) holds. If in addition, we can choose $\widetilde{\beta}_n\to 0$, then for each given $\epsilon$, we have $(1-\widetilde{\beta}_n)\|f-\widehat{f}_{n_1,2}\|_2^2>\frac{1+\alpha_\epsilon}{(1+\alpha'_\epsilon/2)}\|f-\widehat{f}_{n_1,2}\|_2^2$ for some small $\alpha_\epsilon>0$ when $n_1$ is large enough. Now for every $\widetilde{\epsilon}>0$, we can find $\epsilon>0$ such that $3\epsilon\le\widetilde{\epsilon}/3$ and there exists an integer $n_0$ such that when $n\ge n_0$ the probability in (6.9) is upper bounded by $\widetilde{\epsilon}/3$ and $\Delta(\epsilon,n)\le\widetilde{\epsilon}/3$. Consequently when $n\ge n_0$,

$$P\bigg(\frac{\sum_{i=n_1+1}^{n}(f(X_i)-\widehat{f}_{n_1,2}(X_i))^2}{\sum_{i=n_1+1}^{n}(f(X_i)-\widehat{f}_{n_1,1}(X_i))^2}\ge 1+\alpha_\epsilon\bigg)\ge 1-\widetilde{\epsilon}.$$

Recall that we needed the conditions (6.4), (6.5), (6.7) and (6.8) for (6.1) to hold. Under Condition 3, $n_2\|f-\widehat{f}_{n_1,2}\|_2^4/\|f-\widehat{f}_{n_1,1}\|_4^4$ is lower bounded in order in probability by $n_2\|f-\widehat{f}_{n_1,2}\|_2^4/(M_{n_1}^4\|f-\widehat{f}_{n_1,1}\|_2^4)$. From all above, since $\widehat{f}_{n_1,1}$ and $\widehat{f}_{n_1,2}$ converge exactly at rates $p_n$ and $q_n$, respectively, under the $L_2$ loss, we know that for the conclusion of Theorem 1 to hold, it suffices to have these requirements: for every $\epsilon>0$, for some $\widetilde{\beta}_n\to 0$, we have that each of $n_2\widetilde{\beta}_n^2M_{n_1}^{-4}$, $n_2(q_{n_1}/p_{n_1})^4M_{n_1}^{-4}$, $n_2\widetilde{\beta}_nq_{n_1}^2(A_{n_1,\epsilon})^{-2}$, $n_2q_{n_1}^2(A_{n_1,\epsilon})^{-2}$ and $n_2q_{n_1}^2$ goes to infinity.

Under Condition 1, for every $\epsilon>0$, there exists a constant $B_\epsilon>0$ such that $P(\|f-\widehat{f}_{n_1,j}\|_\infty\ge B_\epsilon A_{n_1})\le\epsilon$ when $n_1$ is large enough. That is, for a given $\epsilon>0$, we can take $A_{n_1,\epsilon}=O(A_{n_1})$. Therefore if we have $n_2M_{n_1}^{-4}\to\infty$ and $n_2q_{n_1}^2/(1+A_{n_1})\to\infty$, then we can find $\widetilde{\beta}_n\to 0$ such that the five requirements in the previous paragraph are all satisfied. This completes the proof of Theorem 1.

**Acknowledgments.** Seymour Geisser passed away in March 2004. As is well known, he was one of the pioneers of the widely applicable cross-validation methods. This paper is in memory of his great contributions to statistics and to the School of Statistics at the University of Minnesota. The author thanks a referee, the Associate Editor and the Co-Editor for helpful comments.

School of Statistics
313 Ford Hall
University of Minnesota
224 Church Street S.E.
Minneapolis, Minnesota 55455
USA
E-mail: yyang@stat.umn.edu